# MATCHING NUMBER, HAMILTONIAN GRAPHS AND DISCRETE MAGNETIC LAPLACIANS

JOHN STEWART FABILA-CARRASCO, FERNANDO LLEDÓ, AND OLAF POST

ABSTRACT. In this article, we relate the spectrum of the discrete magnetic Laplacian (DML) on a finite simple graph with two structural properties of the graph: the existence of a perfect matching and the existence of a Hamiltonian cycle of the underlying graph. In particular, we give a family of spectral obstructions parametrised by the magnetic potential for the graph to be matchable (i.e., having a perfect matching) or for the existence of a Hamiltonian cycle. We base our analysis on a special case of the spectral preorder introduced in [FCLP20a] and we use the magnetic potential as a spectral control parameter.

## 1. INTRODUCTION

Spectral graph theory studies the relationship between combinatorial and geometric properties of a graph with the eigenvalues of some matrix associated with it (typically, the adjacency matrix, the combinatorial Laplacian, the signless Laplacian or the normalised Laplacian). Some concrete results in this direction relate the matching number of a graph (i.e., the maximal number of independent edges of a graph) with the spectrum of the combinatorial Laplacian (cf., [MW01]) or the eigenvalues of the signless Laplacian with the circumference of the graph (cf., [WB13]). Moreover, the existence of a Hamiltonian cycle in the graph (i.e., a closed path in a connected graph that contains each vertex exactly once) with bounds in the spectrum of the combinatorial Laplacian (see [Heu95, Moh92]). Recall that deciding if a graph has a Hamiltonian cycle is an NP-complete problem and, therefore, in this context one usually gives sufficient conditions on the spectrum of the Laplacian that guarantee the existence of a Hamiltonian cycle. A different approach is given in [BC10], where the authors show that if the non-trivial eigenvalues of the combinatorial Laplacian are sufficiently close to the average degree, then the graph is Hamiltonian.

The discrete Laplacian can be generalised in a natural way to include a magnetic field which is modelled by a magnetic potential function defined on the set of all directed edges (arcs) of the graph with values in the unit circle, i.e., $\alpha\colon A \to \mathbb{R}/2\pi\mathbb{Z}$. Such operator is called in [Shu94] the Discrete Magnetic Laplacian (DML for short) and is denoted by $\Delta_\alpha$ (see also [Sun94, HS99]). It includes the cases of the combinatorial Laplacian (if $\alpha = 0$) or the signless Laplacian if $\alpha_e = \pi$ on all directed edges $e \in A$ (see Example 2.3). The analysis of the DML is interesting for theoretical aspects and, also, in applications to mathematical physics, particularly in solid state and condensed matter physics, where one uses graphs as a model of a solid (see, e.g., [Est15]).

The analysis of the spectrum of the magnetic Laplacian is a particularly rich object of study because the presence of the magnetic potential amplifies the relationship between the operator and the topology of the graph. In particular, if the graph is connected and has no cycles (e.g., if the graph is a tree), then the magnetic potential has no effect and the DML is unitarily equivalent to the usual combinatorial Laplacian. Besides the evident physical importance of a magnetic field in interaction with a graph, the magnetic potential has many more applications. For example, the magnetic potential $\alpha$ can be interpreted as a Floquet parameter to analyse the Laplacian on an (infinite) periodic graph (see [FCLP18, FCL19, FCLP20a] as well as [KS17, KS19]). Moreover, the magnetic potential plays also the role of a spectral control parameter of the system. In fact, using the magnetic potential as a continuous parameter one can modify the spectrum of the Laplacian and, for instance, raise its first eigenvalue or control the size of spectral gaps in periodic structures (see [FCL19]).

Nevertheless, in the combinatorics literature, the DML is rarely considered since, in principle, the magnetic field is an additional structure of the graph. The aim of this article is to show that the DML with combinatorial weights is useful to address certain questions in discrete mathematics. In particular, we explore the relation between the spectrum of the discrete magnetic Laplacian and two combinatorial properties of the graph: the matching number and the existence of Hamiltonian cycles. We extend some of the results in [Heu95, MW01, WB13] that include statements involving the spectrum of the combinatorial or signless Laplacian. Moreover, the magnetic potential allows us to enlarge the spectral obstructions to the existence of a Hamiltonian cycle in the graph or the existence of a perfect matching.

*Date*: August 11, 2021, 16:21, *File*: `fabila-lledo-post-hamiltonian-ARXIV.tex`.
*Key words and phrases.* Spectral graph theory, discrete magnetic Laplacian, matching number, hamiltonian graph
2020 *Mathematics Subject Classification.* 05C70,05C45,39A70,47A10,05C50 .
JSFC and FLl were supported by Spanish Ministry of Economy and Competitiveness through project DGI MTM2017-84098-P, from the Severo Ochoa Programme for Centres of Excellence in R& D (SEV-2015-0554) and from the Spanish National Research Council, through the *Ayuda extraordinaria a Centros de Excelencia Severo Ochoa* (20205CEX001).

This article is structured as follows: In Section 2, we introduce the notation for the main discrete structures needed. We consider finite and simple graphs (i.e., the graph with a finite number of vertices and with no multiples edges or loops). We include in this section the definition of the DML with combinatorial weights and mention a spectral preorder that controls the spectral spreading of the eigenvalues under edge deletion. We refer to [FCLP20a] for a general analysis of this preorder for multigraphs with general weights and, also, for additional motivation. In Section 3, we introduce some relations between the matching number of the graph and the eigenvalues of the magnetic Laplacian. Moreover, we generalise some spectral upper and lower bounds stated in [Heu95, MW01, WB13] for the combinatorial or signless Laplacian. In Section 4, we address the problem of giving spectral obstructions for the graph being Hamiltonian. In particular, we present examples of graphs where the obstructions given by the usual (or signless) Laplacian in [Heu95, MW01, WB13] do not apply. Nevertheless, for certain values of the of the magnetic potential the DML provides spectral obstructions for the graph to be Hamiltonian.

## 2. Graph theory and spectral preorder

2.1. **Discrete graphs.** A *discrete (oriented) graph* (or, simply, a *graph*) $G = (V, A)$ consists of two disjoint and finite sets $V = V(G)$ and $A = A(G)$, the set of *vertices* and the set of all *oriented edges*, respectively, and a *connection map* $\partial = \partial^G \colon A \longrightarrow V \times V$, where $\partial e = (\partial_- e, \partial_+ e)$ denotes the pair of the *initial* and *terminal* vertex. We also say that $e$ *starts* at $\partial_- e$ and *ends* at $\partial_+ e$. We assume that each oriented edge $e$ (also called *arrow* or *arc*) comes with its oppositely oriented edge $\bar{e}$, i.e., that there is an involution $\bar{\cdot} \colon A \longrightarrow A$ such that $e \neq \bar{e}$ and $\partial_\pm \bar{e} = \partial_\mp e$ for all $e \in A$. If $V(G)$ has $n \in \mathbb{N}$ vertices, we say that $G$ is a *finite* graph of *order* $n$ and we write $|G| = |V(G)| = n$.

We denote by
$$A_v := \{\, e \in A \,|\, \partial_- e = v \,\}$$
the set of all arcs starting at $v$ (alternatively we may also write $A_v^G$). We define the *degree* of the vertex $v$ in the graph $G = (V, A)$ by the cardinality of $A_v$, i.e.,
$$\deg(v) := \deg^G(v) = |A_v|.$$

The *inversion map* gives rise to an action of $\mathbb{Z}_2 = \mathbb{Z}/2\mathbb{Z}$ on the set of arcs $A$. An (unoriented) *edge* is an element of the orbit space $E = A/\mathbb{Z}_2$, i.e., an edge is obtained by identifying the arc $e$ and $\bar{e}$. We denote an (unoriented) graph by $G = (V, E)$ and set $\partial[e] = \{\partial_- e, \partial_+ e\}$ for $[e] \in E$. We say that $[e_1], [e_2] \in E$ *share a vertex* if $\{\partial_- e_1, \partial_+ e_1\} \cap \{\partial_- e_2, \partial_+ e_2\} \neq \emptyset$. To simplify notation, we mostly write $e \in E$ instead of $[e] \in E$ for $e \in A$. Two edges $e_1, e_2 \in E$ in a graph $G = (V, E)$ are *independent* if they are not loops and if they do not share a vertex. A *matching* in a graph is a set of pairwise independent edges. The *matching number* of $G$, denoted by $\mu(G)$, is the cardinality of the maximum number of pairwise independent edges in $G$. A vertex $v$ *belongs to a matching* $M$ if $v \in \bigcup_{e \in M} \partial e$. A *perfect matching* $M$ is a matching where all vertices of $V$ belong to $M$, i.e., where $\bigcup_{e \in M} \partial e = V$. A graph is *matchable* if it has a perfect matching. Recall that if a tree is matchable then the number of its vertices must be even.

2.2. **Magnetic potentials.** Let $G = (V, A)$ be a graph and consider the group $R = \mathbb{R}/2\pi\mathbb{Z}$ with the operation written additively. We consider also the corresponding *cochain groups* of $R$-valued functions on vertices and edges respectively:
$$C^0(G, R) := \{\, \xi \colon V \longrightarrow R \,\} \quad \text{and} \quad C^1(G, R) := \{\, \alpha \colon A \longrightarrow R \,|\, \forall e \in A \colon \alpha_{\bar{e}} = -\alpha_e \,\}.$$
The *coboundary operator* mapping 0-cochains to 1-cochains is given by
$$\mathrm{d} \colon C^0(G, R) \longrightarrow C^1(G, R), \quad \text{where} \quad (\mathrm{d}\xi)_e := \xi(\partial_+ e) - \xi(\partial_- e).$$

**Definition 2.1.** Let $G = (V, A)$ be a graph and $R = \mathbb{R}/2\pi\mathbb{Z}$.
  (i) An *$R$-valued magnetic potential* $\alpha$ is an element of $C^1(G, R)$.
 (ii) We say that $\alpha, \widetilde{\alpha} \in C^1(G, R)$ are *cohomologous* or *gauge-equivalent* and denote this as $\widetilde{\alpha} \sim \alpha$ if $\widetilde{\alpha} - \alpha$ is exact, i.e., if there is $\xi \in C^0(G, R)$ such that $\mathrm{d}\xi = \widetilde{\alpha} - \alpha$, and $\xi$ is called the *gauge*.
(iii) We say that $\alpha$ is *trivial*, if it is cohomologous to 0.

In the sequel, we will omit for simplicity the Abelian group $R$, e.g. we will write $C^1(G)$ instead of $C^1(G, R)$ for the group of magnetic potential etc. We refer to [LP08, Section 5] for additional motivation on homologies of graphs (see also [MY02] and references therein for a version of these homologies twisted by the magnetic potential).

2.3. **The magnetic Laplacian and spectral preorder.** In this section, we will introduce the discrete magnetic Laplacian associated to a graph $G$ with a magnetic potential $\alpha$. We call it simply a *magnetic graph* and denote it by $\mathbf{G} = (G, \alpha)$. We will also introduce a spectral relation between the magnetic Laplacian associated to different graphs.

Given a finite graph $G$, we define the Hilbert space $\ell_2(V) := \{\varphi \colon V \longrightarrow \mathbb{C}\}$ (which is isomorphic to $\mathbb{C}^{|V|}$) and with the inner product defined as usual given by
$$\langle f, g \rangle_{\ell_2(V)} := \sum_{v \in V} f(v) \overline{g(v)}.$$

Note that functions on $V$ may be interpreted as 0-forms (while functions on edges are 1-forms).

**Definition 2.2** (Discrete magnetic Laplacian)**.** Let $G = (V, A)$ be a graph and $\alpha \in C^1(G)$ an $R$-valued magnetic potential, i.e., a map $\alpha \colon A \longrightarrow R$ such that $\alpha_{\bar{e}} = -\alpha_e$ for all $e \in A$, where $R = \mathbb{R}/2\pi\mathbb{Z}$. The *(discrete) magnetic Laplacian* is an operator

$$\Delta_\alpha \colon \ell_2(V) \longrightarrow \ell_2(V) , \qquad (2.1)$$

that acts as

$$(\Delta_\alpha \varphi)(v) = \sum_{e \in A_v} \left( \varphi(v) - \mathrm{e}^{\mathrm{i}\alpha_e} \varphi(\partial_+ e) \right) = \deg(v)\varphi(v) - \sum_{e \in A_v} \mathrm{e}^{\mathrm{i}\alpha_e} \varphi(\partial_+ e). \qquad (2.2)$$

The DML can be seen as a second order discrete operator and one can show that that $\Delta_\alpha$ is positive definite and has spectrum contained in the interval $[0, \max_{v \in V} \deg(v)]$ (see, e.g., [FCLP18, Section 2.3]). If we need to stress the dependence of the DML on the graph $G$, we will write the Laplacian as $\Delta_\alpha^G$. If $G = (V, A)$ is a graph of order $|V(G)| = n$ and magnetic potential $\alpha$, we denote the spectrum of the corresponding magnetic Laplacian by $\sigma(\Delta_\alpha) := \{\lambda_k(\Delta_\alpha) \mid k = 1, \ldots, n\}$. Moreover, we will write the eigenvalues in ascending order and repeated according to their multiplicities, i.e.,

$$0 \leq \lambda_1(\Delta_\alpha) \leq \lambda_2(\Delta_\alpha) \leq \cdots \leq \lambda_n(\Delta_\alpha) .$$

**Example 2.3** (Special cases of the magnetic Laplacian)**.**
  (i) If $\alpha \sim 0$, then $\Delta_\alpha^G$ is unitarily equivalent with the usual combinatorial Laplacian $\Delta_0^G$ (without magnetic potential).
  (ii) Choosing $\alpha_e = \pi$ for all $e \in A$, then $\Delta_\alpha^G = \Delta_\pi^G$ is the *signless* Laplacian.
  (iii) If we choose $R = \{0, \pi\}$, then the magnetic potential is also called *signature*, and $(G, \alpha)$ is called a *signed graph* (see, e.g. [LLPP15] and references therein).

We present a spectral preorder for graphs with magnetic potential which is a particular case the preorder studied for general weights and multigraphs in [FCLP20a]. We refer to this article for additional references, motivation and applications.

**Definition 2.4** (Spectral preorder of magnetic graphs)**.** Let $\mathbf{G} = (G, \alpha)$ and $\mathbf{G}' = (G', \alpha')$ be two magnetic graphs with $|G| = |G'| = n$. We say that $\mathbf{G}$ *is (spectrally) smaller than* $\mathbf{G}'$ *with shift* $r$ and denote this by $\mathbf{G} \preccurlyeq^r \mathbf{G}'$) if

$$\lambda_k(\Delta_\alpha^G) \leq \lambda_{k+r}(\Delta_{\alpha'}^{G'}) \qquad \text{for all } 1 \leq k \leq n - r.$$

If $r = 0$ we write again simply $\mathbf{G} \preccurlyeq \mathbf{G}'$.

The relation $\preccurlyeq$ is a preorder (i.e., a reflexive and transitive relation) on the class of magnetic graphs (cf. [FCLP20a, Proposition 3.11]). Moreover, $\mathbf{G}' \preccurlyeq^s \mathbf{G} \preccurlyeq^r \mathbf{G}'$ means that $\mathbf{G}' \preccurlyeq^s \mathbf{G}$ and $\mathbf{G} \preccurlyeq^r \mathbf{G}'$. In particular, if $s = 0$ and $r = 1$, the relation $\mathbf{G}' \preccurlyeq \mathbf{G} \preccurlyeq^1 \mathbf{G}'$ describes the usual interlacing of eigenvalues:

$$\lambda_1(\Delta_{\alpha'}^{G'}) \leq \lambda_1(\Delta_\alpha^G) \leq \lambda_2(\Delta_{\alpha'}^{G'}) \leq \lambda_2(\Delta_\alpha^G) \leq \cdots \leq \lambda_{n-1}(\Delta_{\alpha'}^{G'}) \leq \lambda_{n-1}(\Delta_\alpha^G) \leq \lambda_n(\Delta_{\alpha'}^{G'}).$$

In the following result we apply the spectral preorder to control the spectral spreading of the DML due to edge deletion and keeping the same magnetic potential on the remaining edges. Recall that for a graph $G = (V, E)$ and an edge $e_0 \in E$, we denote by $G - e_0$ the graph given by $G - e_0 = (V, E \setminus \{e_0\})$.

**Theorem 2.5.** *Let $\mathbf{G} = (G, \alpha)$ and $\mathbf{G}' = (G', \alpha')$ be two magnetic graphs where $G' = G - e_0$ for some $e_0 \in E(G)$, and $\alpha_e = \alpha'_e$ for all $e \in A(G')$, then*

$$\mathbf{G}' \preccurlyeq \mathbf{G} \preccurlyeq^1 \mathbf{G}' .$$

Applying several times the previous relations it is clear that if $G'$ is obtained from $G$ by deleting $r$ edges, i.e., if $G' = G - \{e_1, \ldots e_r\}$, then

$$\mathbf{G}' \preccurlyeq \mathbf{G} \preccurlyeq^r \mathbf{G}' . \qquad (2.3)$$

The preceding theorem generalises to the DML with arbitrary magnetic potential some known interlacing results, namely [Heu95, Lemma 2] (for combinatorial and signless Laplacians) or [Moh91, Theorem 3.2] and [Fie73, Corollary 3.2] (for combinatorial Laplacian). The general case for any magnetic potential and arbitrary weights is proved in [FCLP20a, Theorem 4.1], for the normalised Laplacian in [But07, Theorem 1.1] and for the standard Laplacian in [CDH04, Theorem 2.2].

## 3. Matching number and the discrete magnetic Laplacian

In this section, we relate some bounds for the eigenvalues of the magnetic Laplacian with the matching number of the underlying graph. In particular, we give a spectral obstruction provided by the DML to the existence of a perfect matching of the graph. In certain examples this obstruction is not effective for the combinatorial and signless Laplacian, but it works for the DML with a certain non-trivial magnetic potential. In this sense the presence of the magnetic potential (thought of as a continuous parameter) makes the spectral obstruction applicable for many more cases (see Example 3.7).

We begin considering the case of trees, so that any DML is unitarily equivalent to usual combinatorial Laplacian. The first result essentially says that if $T$ is a tree with matching number $\mu(T)$, then the first $\mu(T)+1$ eigenvalues are smaller or equal than 2 and at least $\mu(T)$ eigenvalues are greater or equal than 2. Recall that we write the eigenvalues in ascending order and repeated according to their multiplicities.

**Theorem 3.1.** *Let $T$ be a tree on $n$ vertices and matching number $\mu(T)$, then*
$$\lambda_{\mu(T)+1}(\Delta^T) \leq 2 \leq \lambda_{n-\mu(T)+1}(\Delta^T) \,.$$

*Proof.* Consider the graph $G$ as the disjoint union of $\mu(T)$ complete graphs $K_2$ and $n - 2\mu(T)$ many isolated vertices, i.e.,
$$G = \left( \bigcup_{i=1}^{\mu(T)} K_2 \right) \cup \left( \bigcup_{i=1}^{n-2\mu(T)} K_1 \right) \,.$$

Then $G$ is a graph on $n$ vertices with $\sigma(\Delta^G) = \{0^{(n-\mu(T))}, 2^{(\mu(T))}\}$, where the superscripts denote the multiplicities of each eigenvalue. In particular,
$$\lambda_{n-\mu(T)+k}(\Delta^G) = 2 \quad \text{for all} \quad k \in \{1,2,\ldots,\mu(T)\} \,. \tag{3.1}$$

The graph $G$ has $\mu(T)$ edges and is obtained from the tree $T$ (which has $n-1$ edges) by deleting the $(n-1-\mu(T))$ edges that do not belong to the matching. Then, by Theorem 2.5, we obtain
$$G \preccurlyeq T \stackrel{n-\mu(T)-1}{\preccurlyeq} G \,. \tag{3.2}$$

For $k=1$ in Eq. (3.1) together with the left relation of Eq. (3.2), it follows that $2 = \lambda_{n-\mu(T)+1}(\Delta^G) \leq \lambda_{n-\mu(T)+1}(\Delta^T)$. Similarly, from the right relation of Eq. (3.2) applied to the case $k = \mu(T)$ in Eq. (3.1) we obtain $\lambda_{\mu(T)+1}(\Delta^T) \leq \lambda_n(\Delta^G) = 2$. □

We mention next some easy consequences of the preceding theorem. Recall first that if $T$ is matchable, then one and only one eigenvalue is equal to 2 (see [MW01, Theorem 2]). This fact follows immediately from the preceding theorem.

**Corollary 3.2.** *Let $T$ be a matchable tree on $n$ vertices, then $n$ is even and*
$$\lambda_{\frac{n}{2}+1}(\Delta^T) = 2 \,.$$

*Proof.* If $T$ is a matchable tree on $n$ vertices, then $n$ is an even number and $\mu(T) = \frac{n}{2}$. By Theorem 3.1 we obtain that $\lambda_{\frac{n}{2}+1}(\Delta^T) \leq 2 \leq \lambda_{\frac{n}{2}+1}(\Delta^T)$ and hence $\lambda_{\frac{n}{2}+1}(\Delta^T) = 2$. □

**Corollary 3.3.** *Let $T$ be a matchable tree on $n$ vertices, then $n$ is even and*
$$\lambda_{\frac{n}{2}}(\Delta^T) < 2 < \lambda_{\frac{n}{2}+2}(\Delta^T) \,.$$

*Proof.* Because $T$ is a matchable tree, by Corollary 3.2, it follows that $\lambda_{\frac{n}{2}+1}(\Delta^T) = 2$ is an eigenvalue. In general, if $T$ is a tree, we know from [GRS90, Theorem 2.1(ii)] that the multiplicity of any integer eigenvalue larger that 1 is exactly equal to one and, therefore, the result follows. □

We use next the spectral preorder to show the following spectral bounds. Note that the second inequality is already shown in [MW01].

**Theorem 3.4.** *Let $T$ be a tree on $n$ vertices. If $n > 2\mu(T)$, then*
$$\lambda_{\mu(T)+1}(\Delta^T) < 2 < \lambda_{n-\mu(T)+1}(\Delta^T) \,.$$

*Proof.* From Theorem 3.1 we already have $\lambda_{\mu(T)+1}(\Delta^T) \leq 2 \leq \lambda_{n-\mu(T)+1}(\Delta^T)$.

Suppose that $n$ is an odd number. If the value 2 is an eigenvalue, then it is an integer eigenvalue. From [GRS90, Theorem 2.1(i)] we conclude that 2 divides $n$ giving a contradiction. Therefore, 2 is not an eigenvalue and follows that $\lambda_{\mu(T)+1}(\Delta^T) < 2 < \lambda_{n-\mu(T)+1}(\Delta^T)$.

Suppose that $n$ is an even number. If $n > 2\mu(T)$, then there exists a vertex $v_0 \in V(T)$ such that $v_0$ does not belong to the maximum matching. Let $r = \deg(v_0)$ and denote by $T - v_0$ the graph obtained from $T$ by deleting the vertex $v_0$ and all $r$ edges adjacent to $v_0$. Denote the connected components of $T - v_0$ by $T_1, T_2, \cdots, T_r$. We can assume that $T_1$ has an odd number of vertices, since otherwise $n = \sum_{i=1}^{r} |V(T_i)| + 1$ is odd giving a contradiction. Hence, the graph $T - V(T_1)$ has also an odd number of vertices. Define $n_1 = |T_1|$, which is an odd number and by the previous case it follows
$$\lambda_{\mu(T_1)+1}(\Delta^{T_1}) < 2 < \lambda_{n_1-\mu(T_1)+1}(\Delta^{T_1}) \,.$$

Similarly, $|T - V(T_1)| = n - n_1$ is an odd number and, again,
$$\lambda_{\mu(T-V(T_1))+1}(\Delta^{T-V(T_1)}) < 2 < \lambda_{n-n_1-\mu(T-V(T_1))+1}(\Delta^{T-V(T_1)}) \,.$$

From Theorem 2.5 we obtain for the disjoint union of graphs

$$T_1 \cup (T - V(T_1)) \preccurlyeq T \overset{1}{\preccurlyeq} T_1 \cup (T - V(T_1))$$

as $T_1 \cup (T - V(T_1))$ is obtained from $T$ by deleting one edge. Moreover, we have $\mu(T) = \mu(T_1) + \mu(T - V(T_1))$ and we conclude $\lambda_{\mu(T)+1}(\Delta^T) < 2 < \lambda_{n-\mu(T)+1}(\Delta^T)$. □

We focus next on general finite simple graphs $G$. It is a well known fact that for any matching $M$ of $G$ there exists a spanning tree $T$ of $G$ which includes all the edges of $M$.

**Corollary 3.5.** *Let $G$ be a connected graph on $n$ vertices, $m$ edges and matching number $\mu(G)$. For any magnetic potential $\alpha$ we have*

$$\lambda_{\mu(G)+n-m}(\Delta_\alpha^G) \leq 2 \leq \lambda_{n-\mu(G)+1}(\Delta_\alpha^G) \;.$$

*Proof.* Consider a spanning tree $T$ of $G$ with the same matching number than $G$, i.e., $\mu(T) = \mu(G)$. Then, $T$ is a graph (with $n-1$ edges) obtained from $G$ by deleting $m-(n-1)$ edges and by Theorem 2.5 we obtain

$$T \preccurlyeq \mathbf{G} \overset{m-n+1}{\preccurlyeq} T \;.$$

By Theorem 3.1 together with the previous relation and the fact that all magnetic Laplacians on a tree are unitarily equivalent (see Example 2.3 (i)) we have

$$2 \leq \lambda_{n-\mu(T)+1}(\Delta^T) \leq \lambda_{n-\mu(G)+1}(\Delta_\alpha^G) \quad \text{and} \quad \lambda_{\mu(G)+n-m}(\Delta_\alpha^G) \leq \lambda_{\mu(T)+1}(\Delta^T) \leq 2$$

for any magnetic potential $\alpha$ on $G$ concluding the proof. □

The next corollary gives a simple family of spectral obstructions for the graph being matchable.

**Corollary 3.6.** *Let $G$ be a graph on $n$ vertices where $n$ is even. If there exists a magnetic potential $\alpha$ such that*

$$\lambda_{\frac{n}{2}+1}(\Delta_\alpha^G) < 2 \;,$$

*then $G$ is not-matchable.*

*Proof.* Suppose that $G$ is a matchable graph, then $\mu(G) = \frac{n}{2}$ and from Corollary 3.5 we conclude

$$2 \leq \lambda_{n-\mu(G)+1}(\Delta_\alpha^G) = \lambda_{\frac{n}{2}+1}(\Delta_\alpha^G)$$

which gives a contradiction. Therefore $G$ is not-matchable. □

**Example 3.7.** Consider the graph $G$ given in Figure 1. The spectrum of the combinatorial Laplacian (i.e., with magnetic potential $\alpha \sim 0$) is

$$\sigma(\Delta_0^G) = \left\{0, 3-\sqrt{5}, 1, 2, 3, \sqrt{5}+3\right\} \;.$$

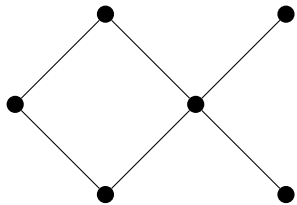

FIGURE 1. Example of a bipartite graph that is not matchable. The spectral obstruction holds only for non-trivial values of the magnetic potential.

Observe that the graph $G$ is bipartite, so that the signless combinatorial Laplacian is unitarily equivalent to the usual combinatorial Laplacian. Therefore, if $\alpha_e = \pi$ for each edge, then the spectra of the combinatorial and signless combinatorial Laplacians coincide, i.e., $\sigma(\Delta_0^G) = \sigma(\Delta_\pi^G)$. In particular, $\lambda_{\frac{n}{2}+1}(\Delta_0^G) = \lambda_{\frac{n}{2}+1}(\Delta_\pi^G) = \lambda_4(\Delta_0^G) = 2$ and, therefore, the eigenvalues of the combinatorial and signless Laplacians provide no obstruction to the matchability of $G$. But, if we consider the magnetic potential $\alpha'$ with value equal to $\pi$ only on one edge of the cycle and zero everywhere else, then the spectrum is given by

$$\sigma(\Delta_{\alpha'}^G) \approx \{0.23, 0.58, 1, 1.63, 3.41, 5.12\} \;.$$

In particular, $\lambda_4(\Delta_{\alpha'}^G) < 2$ and therefore we can conclude from Corollary 3.6 that $G$ is not matchable. In this example any non-trivial magnetic potential $\alpha \nsim 0$ provides a spectral obstruction since $\lambda_4(\Delta_\alpha^G) < 2$ as the following plot of the eigenvalue $\lambda_4(\Delta_\alpha^G)$ for different values of $\alpha$ shows.

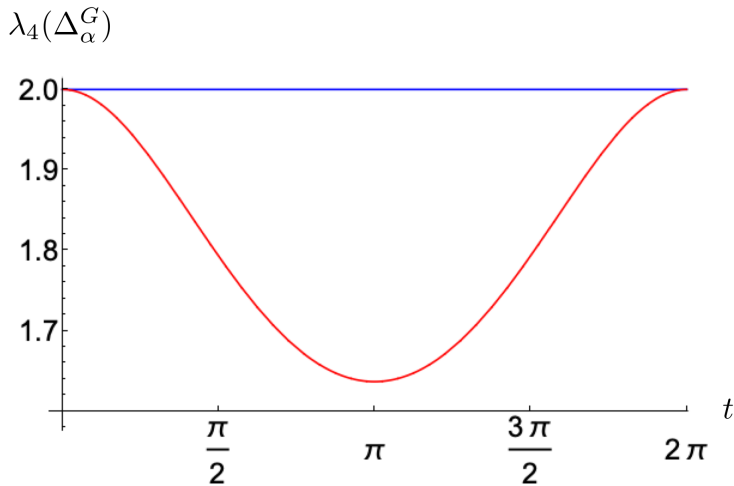

FIGURE 2. Plot of the eigenvalue $\lambda_4(\Delta_\alpha^G)$ as a function of a magnetic potential $\alpha_e = t$ for one edge of the cycle in the graph in Figure 1 and zero everywhere else.

The next proposition generalises to the DML with arbitrary magnetic potential $\alpha$ results known for the combinatorial and signless Laplacians (see [MW01, Theorem 4] and [WB13, Lemma 2.4]).

**Proposition 3.8.** *Let $G$ a connected graph with $n$ vertices, $m$ edges. Moreover, let $\alpha$ be a magnetic potential on $G$.*
  (i) *If $n > 2\mu(G)$, then $\lambda_{\mu(G)+n-m}(\Delta_\alpha^G) < 2 < \lambda_{n-\mu(G)+1}(\Delta_\alpha^G)$.*
  (ii) *If $n = 2\mu(G)$, then $\lambda_{\frac{3n}{2}-m-1}(\Delta_\alpha^G) < 2 < \lambda_{\frac{n}{2}+2}(\Delta_\alpha^G)$.*

*Proof.* Consider a spanning tree $T$ of $G$ with the same matching number than $G$, i.e., $\mu(T) = \mu(G)$. Then, $T$ is a graph (with $n-1$ edges) obtained from $G$ by deleting $m - (n-1)$ edges and by Theorem 2.5 we obtain

$$T \preccurlyeq \mathbf{G} \stackrel{m-n+1}{\preccurlyeq} T. \tag{3.3}$$

First note that the preceding relations in Eq. (3.3) together with Theorem 3.4 give

$$\lambda_{\mu(G)+n-m}(\Delta_\alpha^G) \leq \lambda_{\mu(T)+1}(\Delta^T) < 2 < \lambda_{n-\mu(T)+1}(\Delta^T) \leq \lambda_{n-\mu(G)+1}(\Delta_\alpha^G).$$

Second, Corollary 3.3 together with Theorem 3.4 gives

$$\lambda_{\frac{3n}{2}-m-1}(\Delta_\alpha^G) \leq \lambda_{\frac{n}{2}}(\Delta^T) < 2 < \lambda_{\frac{n}{2}+2}(\Delta^T) \leq \lambda_{n-\mu(G)+2}(\Delta_\alpha^G)$$

which concludes the proof. □

## 4. Hamiltonian graphs and the magnetic potential

A cycle which contains every vertex of the graph is called a *Hamiltonian cycle* and a graph is said to be *Hamiltonian* if it has a Hamiltonian cycle. Some results that connect the existence of a Hamilton cycle in the graph and bounds of the eigenvalues of the combinatorial Laplacian are given in [Heu95, Moh92]. In particular, the next result generalise the main result of Theorems 1 and 1' in [Heu95].

**Theorem 4.1.** *Let $\mathbf{G} = (G, \alpha)$ be a magnetic graph with $n$ vertices, $m$ edges and with a magnetic potential $\alpha$. If $G$ contains a Hamiltonian cycle $C_n$, then*

$$\mathbf{C}_n \preccurlyeq \mathbf{G} \stackrel{m-n}{\preccurlyeq} \mathbf{C}_n$$

*where $\mathbf{C}_n = (C_n, \alpha')$ is the magnetic graph with $\alpha'$ denoting the restriction of $\alpha$ to the edges of the cycle, i.e., $\alpha' = \alpha\!\restriction_{E(C_n)}$.*

*Proof.* Since $C_n$ is obtained from $G$ by deleting $m - n$ edges we obtain applying Eq. (2.3) the required inequalities (see also Theorem 2.5). □

As a consequence of the previous result we mention a first spectral obstruction of the DML to the existence of a Hamiltonian cycle.

**Corollary 4.2.** *Let $G = (V, A)$ be a graph on $n$ vertices. Assume that there exists an index $k \in \{1, 2, \ldots, n\}$ and a constant magnetic potential $\alpha = t$ (i.e., there is $t \in \mathbb{R}/2\pi\mathbb{Z}$ with $\alpha_e = t$ for all $e \in A$) such that*
$$\lambda_k(\Delta_{\alpha'}^{C_n}) > \lambda_k(\Delta_{\alpha}^{G}),$$
*where $C_n$ is the cycle on $n$ vertices and $\alpha' = t$ is the magnetic potential on $C_n$. Then $G$ is non-Hamiltonian.*

*Proof.* The inequality follows directly from the relation $C_n \preccurlyeq G$ of Theorem 4.1. □

**Example 4.3.** Consider the graph $G$ on 6 vertices given in Figure 3 with trivial magnetic potential $\alpha = 0$, i.e., $\mathbf{G} = (G, 0)$. Then the spectrum of the combinatorial Laplacian is given by:
$$\sigma(\Delta_0^G) = \left\{0, \frac{1}{2}\left(7 - \sqrt{17}\right), 2, 3, 4, \frac{1}{2}\left(\sqrt{17} + 7\right)\right\} \approx \{0, 1.43, 2, 3, 4, 5.56\}.$$
The spectrum of the cycle on 6 vertices $C_6$ with magnetic potential $\alpha \sim 0$ is
$$\sigma(\Delta_0^{C_6}) = \{0, 1, 1, 3, 3, 4\}.$$

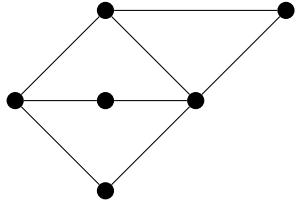

FIGURE 3. Example of a graph that can be shown to be non-Hamiltonian using the eigenvalues of the magnetic Laplacian.

It follows that $\mathbf{C}_6 \preccurlyeq \mathbf{G}$ and the spectra of the combinatorial Laplacian gives no obstruction to the existence of a Hamiltonian cycle of $G$. Similarly, the analysis for the signless Laplacian gives no obstruction. In fact, consider the graph with magnetic potential $\alpha_e = \pi$ for all $e \in A$, i.e., $\mathbf{G} = (G, \alpha = \pi)$. Then, the eigenvalues of the signless Laplacian are
$$\sigma(\Delta_\pi^G) \approx \{0.30, 1.22, 2, 3, 3.58, 5.87\},$$
while the spectrum of the $C_6$ for the signless Laplacian coincides with the spectrum of the usual Laplacian since $C_6$ is bipartite, i.e.,
$$\sigma(\Delta_\pi^{C_6}) = \{0, 1, 1, 3, 3, 4\}.$$
Again, it follows that $\mathbf{C}_6 \preccurlyeq \mathbf{G}$ and the spectrum of the signless Laplacians gives no obstruction.

Consider now the constant magnetic potential $\alpha_e = \pi/2$ for all $e \in A$. The spectrum of the corresponding DML provides the obstruction. In fact, the spectrum of the magnetic Laplacian associated to the magnetic graph $\mathbf{G} = (G, \alpha = \pi/2)$ is given by
$$\sigma(\Delta_{\pi/2}^G) \approx \{0.13, 1.35, 2, 3, 3.77, 5.73\}.$$
and the spectrum of the cycle $C_6$ with constant magnetic potential $\alpha_e = \pi/2$, $e \in A$, is
$$\sigma(\Delta_{\pi/2}^{C_6}) \approx \{0.25, 0.52, 2.39, 2.83, 3.98\}.$$
It is clear that $\lambda_1(\Delta_{\pi/2}^{C_6}) \approx 0.25 > 0.13 = \lambda_1(\Delta_{\pi/2}^G)$, hence by Corollary 4.2 we conclude that $G$ is non-Hamiltonian. Note that there are other values of the magnetic potential where the first eigenvalue provides a similar obstruction (see Figure 4).

As a consequence of Theorem 4.1, we generalise next a result for the signless Laplacian proved in [WB13, Theorem 2.8] to arbitrary DMLs.

**Corollary 4.4.** *Let $G$ be a Hamiltonian graph of order $n > 3$, and $\alpha$ any magnetic potential on $G$. Then,*
  (i) *if $n$ is even, then $2 \leq \lambda_{\frac{n}{2}+1}(\Delta_\alpha^G)$ and $2 < \lambda_{\frac{n}{2}+2}(\Delta_\alpha^G)$.*
  (ii) *if $n$ is odd, then $\lambda_{\frac{n+1}{2}+1}(\Delta_\alpha^G) > 2$*

*Proof.* (i) If $n$ be even then $G$ is matchable and $\mu(G) = \frac{n}{2}$. By Corollary 3.5 we conclude $2 \leq \lambda_{n-\mu(G)+1}(\Delta_\alpha^G) = \lambda_{\frac{n}{2}+1}(\Delta_\alpha^G)$ and by Proposition 3.8 it follows that $\lambda_{\frac{n}{2}+2}(\Delta_\alpha^G) > 2$.
  (ii) If $n$ is odd then $\mu(G) = \frac{n-1}{2}$ and by Proposition 3.8 we conclude $2 < \lambda_{n-\frac{n-1}{2}+1}(\Delta_\alpha^G)$. □

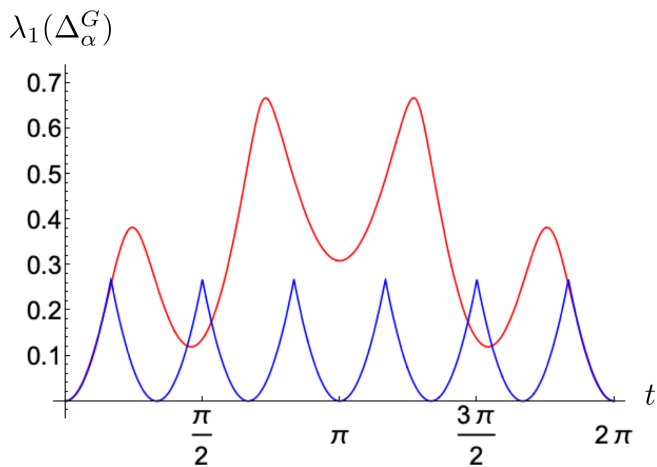

FIGURE 4. First eigenvalue $\lambda_1$ for the Laplacians $\Delta_\alpha^G$ (red) and $\Delta_\alpha^{C_6}$ (blue) for different values of the magnetic potential $\alpha_e = t$, $e \in A$. There is an open interval of values around $t = \pi/2$ and $t = 3\pi/2$ that provide spectral obstructions of the DML to the existence of Hamiltonian cycles.

The previous result gives a sufficient condition for a graph $G$ to be Hamiltonian. In particular if $G$ is a graph with even vertices and $\alpha$ a magnetic potential such that
$$\lambda_{\frac{n}{2}+1}(\Delta_\alpha^G) < 2$$
then $G$ is non-Hamiltonian. We apply this reasoning in the following example.

**Example 4.5.** Consider the graph given in Figure 1 and recall that in Example 3.7, we showed that for any $\alpha \not\sim 0$
$$\lambda_4(\Delta_\alpha^G) < 2 \tag{4.1}$$
concluding that the graph is not matchable. Now for this example we conclude from Corollary 4.4 (i) that one can use again inequality Eq. (4.1) to show that $G$ is non-Hamiltonian either.

Department of Mathematics, University Carlos III de Madrid, Avda. de la Universidad 30, 28911. Leganés (Madrid), Spain and Instituto de Ciencias Matemáticas (CSIC-UAM-UC3M-UCM), Madrid
*Email address*: jfabila@math.uc3m.es

Department of Mathematics, University Carlos III de Madrid, Avda. de la Universidad 30, 28911. Leganés (Madrid), Spain and Instituto de Ciencias Matemáticas (CSIC-UAM-UC3M-UCM), Madrid
*Email address*: flledo@math.uc3m.es

Fachbereich 4 – Mathematik, Universität Trier, 54286 Trier, Germany
*Email address*: olaf.post@uni-trier.de